\documentclass[12pt]{amsart}
\usepackage[]{fontenc}
\usepackage[latin1]{inputenc}
\usepackage{amssymb}

\makeatletter

 \theoremstyle{plain}
 \theoremstyle{definition}
  \newtheorem{example}{Example}

\date{}
\usepackage{mathrsfs}
\usepackage{hyperref}

\makeatother
\begin{document}

\newcommand{\set}[2]{\left\{  #1\, |\, #2\right\}  }

\newcommand{\cyc}[1]{\mathbb{Q}\left[\zeta_{#1}\right]}

\newcommand{\Mod}[3]{#1\equiv#2\, \left(\mathrm{mod}\, \, #3\right)}

\newcommand{\FA}[1]{\mathcal{X}_{2}}

\newcommand{\FB}[1]{\mathcal{X}_{1}}

\newcommand{\tr}{\mathrm{Tr}}

\newcommand{\xrep}{\mathscr{M}\mathcal{\left(X\right)}}

\newcommand{\cl}{\mathscr{C}\ell\left(\mathcal{X}\right)}

\newcommand{\rep}{\mathrm{Xrep}}

\newcommand{\FC}[1]{}

\newcommand{\FD}[1]{}

\newcommand{\FE}[2]{}

\newcommand{\id}{\mathrm{\mathbf{id}}}

\newcommand{\FF}[2]{}

\newcommand{\irr}{\mathrm{Irr}\left(\mathcal{X}\right)}

\newcommand{\triv}{\mathbf{1}}

\newcommand{\vac}{\mathbf{\underline{0}}}

\newcommand{\idr}{\mathbf{\underline{1}}}

\title{Characters of crossed modules and premodular categories}

\author{P. Bantay}

\address{Institute for Theoretical Physics, E\"otv\"os  University, Budapest}

\begin{abstract}
A general procedure is presented which associates to a finite crossed module a premodular category, generalizing the representation categories of  a finite group and of its double, and the extent to which the resulting  category fails to be modular is explained.
\end{abstract}
\maketitle

\section{Introduction}

Modular Tensor Categories \cite{Turaev,Bakalov-Kirillov} have attracted
much attention in recent years, which is due to the recognition of
their importance in both pure mathematics (3-dimensional topology,
representations of Vertex Operator Algebras ) and  theoretical physics
(Rational Conformal Field Theory, Topological Field Theories). They
are also closely related to Moonshine \cite{Convay-Norton,Borcherds2,25years}:
a most interesting (and mysterious) example of a Modular Tensor Category,
which is responsible for some of the deeper aspects of Moonshine,
is the MTC associated to the Moonshine orbifold, i.e. the fixed point
VOA of the Moonshine module under the action of the Monster.

As in every branch of science, a deeper understanding of Modular Tensor
Categories requires a suitable supply of examples. Since the work
of Huang \cite{Huang}, we know that the module category of any rational
VOA (satisfying some technical conditions) is modular, but this important
result doesn't help us that much, because VOA-s are pretty complicated
objects usually hard to deal with. This leads to the desire of associating
MTC-s to simpler and more accessible algebraic objects. There are
several such constructions, a most notable being the one that associates
to a finite group the module category of its (Drinfeld) double \cite{DPR,Hopf1}.
The aim of the present note is to sketch a generalization of this
last construction, associating to any (finite) crossed module a premodular
category, i.e. a braided tensor category that falls short of being
modular. The idea behind is to use 'higher dimensional groups', whose
simplest instance are crossed modules \cite{Wh1,Brown1}, for constructing
Modular Tensor Categories. In the sequel we'll examine to which extent
this idea may be put to work.

The plan of the paper is the following. In the next section we'll
recall some basic definitions and results about crossed modules. In
Section \ref{sec:The-category} we introduce our basic object of study,
the tensor category associated to the crossed module, and discuss
some of its properties. Section \ref{sec:Characters} describes the
notion of characters of crossed modules, the main technical tool in
our study. Section \ref{sec:premodular} discusses the premodular
structure of the category, and the extent to which it fails to be
modular. We conclude by some remarks on the possible applications
of the results presented.

We have decided to present only an outline of the theory, without
going into detailed proofs, since we felt that their inclusion wouldn't
help to clarify the arguments, but could hide the main line of thought.
Detailed proofs of all the results to be presented could be supplied
by exploiting the close analogy with the character theory of finite
groups.

\section{Crossed modules\label{sec:Crossed-modules}}

To begin with, let's recall that an action of the group $G$ on the
group $M$ is a homomorphism $G\rightarrow\mathrm{Aut}\left(M\right)$
or, what is the same, a map $\mu:M\times G\rightarrow M$ such that

\begin{enumerate}
\item $\mu\left(m_{1}m_{2},g\right)=\mu\left(m_{1},g\right)\mu\left(m_{2},g\right)$
for all $m_{1},m_{2}\in M$ and $g\in G$;
\item $\mu\left(m,g_{1}g_{2}\right)=\mu\left(\mu\left(m,g_{1}\right),g_{2}\right)$
for all $m\in M$ and $g_{1},g_{2}\in G$.
\end{enumerate}
As is customary, we'll use the exponential notation $\mu\left(m,g\right)=m^{g}$
in the sequel.

A crossed module \cite{Wh1,Brown1,Gilbert} is nothing but a 4-tuple
$\mathcal{X}=\left(\FB{},\FA{},\mu,\partial\right)$, where $\FB{},\FA{}$
are groups, $\mu$ is an action of $\FB{}$ on $\FA{}$, and $\partial:\FA{}\rightarrow\FB{}$
is a homomorphism, called the boundary map, that satisfies

\begin{description}
\item [XMod1\label{xmod1}]: $\partial\left(m^{g}\right)=g^{-1}\left(\partial m\right)g$
for all $m\in\FA{}$ and $g\in\FB{}$;
\item [XMod2\label{xmod2}]: $m^{\partial n}=n^{-1}mn$ for all $m,n\in\FA{}$.
\end{description}
A crossed module is finite if both $\FB{}$ and $\FA{}$ are finite
groups. Examples of crossed modules abound in algebra and topology,
let's just cite two, coming from group theory, that will guide our
investigations later.

\begin{example}
\label{ex:RG}For a group $G$, we'll denote by $\mathcal{R}G$ the
crossed module $\left(G,\triv,\mu,\partial\right)$, where $\triv$
denotes the trivial subgroup of $G$, i.e. $\triv=\{1\}$, and both
the action $\mu$ and the boundary map $\partial$ are trivial.
\end{example}

\begin{example}
\label{ex:DG}If $G$ is a group, $\mathcal{D}G$ is the crossed module
$\left(G,G,\mu,\id\right)$, where $\mu$ is the conjugation action,
i.e. $\mu\left(m,g\right)=g^{-1}mg$, and $\id:g\mapsto g$ is the
trivial map.
\end{example}
A standard consequence of the defining properties of a crossed module
is that $K=\ker\partial$ is a central subgroup of $\FA{}$, $I=\mathrm{\mathrm{im}}\,\partial$
is a normal subgroup of $\FB{}$, and one has an exact sequence \begin{equation}
\triv\rightarrow K\rightarrow\FA{}\rightarrow\FB{}\rightarrow C\rightarrow\triv\label{eq:exseq}\end{equation}
where $C=\FB{}/I$ is the cokernel of $\partial$ \cite{Brown1}.
In particular, $\left|\FA{}\right|\left|C\right|=\left|K\right|\left|\FB{}\right|$
for a finite crossed module.

Finally, a morphism $\phi:\mathcal{X}\rightarrow\mathcal{Y}$ between
the crossed modules $\mathcal{X}=\left(\FB{},\FA{},\mu_{\mathcal{X}},\partial_{\mathcal{X}}\right)$
and $\mathcal{Y}=\left(\mathcal{Y}_{1},\mathcal{Y}_{2},\mu_{\mathcal{Y}},\partial_{\mathcal{Y}}\right)$
is a pair $\left(\phi_{1},\phi_{2}\right)$, where $\phi_{i}:\mathcal{X}_{i}\rightarrow\mathcal{Y}_{i}$
are group homomorphisms for $i=1,2$, and the following relations
hold:\begin{eqnarray*}
\partial_{\mathcal{Y}}\circ\phi_{1} & = & \phi_{2}\circ\partial_{\mathcal{X}}\\
\mu_{\mathcal{Y}}\circ\left(\phi_{2}\times\phi_{1}\right) & = & \phi_{1}\circ\mu_{\mathcal{X}}\end{eqnarray*}

\section{The category\label{sec:The-category}}

To any finite crossed module $\mathcal{X}=\left(\FB{},\FA{},\mu,\partial\right)$
we'll associate a braided tensor category $\xrep$, which falls short
of being modular. Let's begin by describing the objects and morphisms
of $\xrep$. Here and in the sequel, we use the notation \[
\delta\left(x,y\right)=\begin{cases}
1 & \mathrm{if}\,\, x=y,\\
0 & \mathrm{otherwise.}\end{cases}\]

An object of $\xrep$ is a triple $\left(V,P,Q\right)$, where $V$
is a complex linear space, while $P$ and $Q$ are maps $P:\FA{}\rightarrow\mathrm{End}\left(V\right)$
and $Q:\FB{}\rightarrow\mathrm{GL}\left(V\right)$ such that \begin{eqnarray}
P\left(m\right)P\left(n\right) & = & \delta\left(m,n\right)P\left(m\right)\label{eq:Pidemp}\\
\sum_{m\in\FA{}}P\left(m\right) & = & \id_{V}\label{eq:Psum}\\
Q\left(g\right)Q\left(h\right) & = & Q\left(gh\right)\label{eq:Qprod}\\
P\left(m\right)Q\left(g\right) & = & Q\left(g\right)P\left(m^{g}\right)\label{eq:PQcommutator}\end{eqnarray}
By the dimension of an object $\left(V,P,Q\right)$ we'll mean the
dimension of the linear space $V$. A morphism $\phi:\left(V_{1},P_{1},Q_{1}\right)\rightarrow\left(V_{2},P_{2},Q_{2}\right)$
between two objects of $\xrep$ is a linear map $\phi:V_{1}\rightarrow V_{2}$
such that $\phi\circ P_{1}\left(m\right)=P_{2}\left(m\right)\circ\phi$
for all $m\in\FA{}$ and, $\phi\circ Q_{1}\left(g\right)=Q_{2}\left(g\right)\circ\phi$
for all $g\in\FB{}$. In general, we won't distinguish isomorphic
objects of $\xrep$.

Let's look at a couple of illustrating examples of objects of $\xrep$
for a finite crossed module $\mathcal{X}=\left(\FB{},\FA{},\mu,\partial\right)$.

\begin{example}
\label{ex:identity}The triple $\idr=\left(V,P,Q\right)$, with $V=\mathbb{C}$,
$P\left(m\right)=\delta\left(m,1\right)\id_{V}$ and $Q\left(g\right)=\id_{V}$,
is a one dimensional object of $\xrep$, that we'll call the trivial
object.
\end{example}

\begin{example}
The triple $\mathbf{\underline{R}}=\left(V,P,Q\right)$, with $V=\mathbb{C}\left(\FB{}\times\FA{}\right)$
and $P\left(m\right)\phi:\left(x,y\right)\mapsto\delta\left(m,y^{x}\right)\phi\left(x,y\right)$,
$Q\left(g\right)\phi:\left(x,y\right)\mapsto\phi\left(xg,y\right)$
for $\phi\in V$, is an object of $\xrep$, that we'll call the regular
object. Clearly, $\dim\mathbf{\underline{R}}=\left|\FB{}\right|\left|\FA{}\right|$.
\end{example}

\begin{example}
\label{ex:vacuum}The triple $\vac=\left(V,P,Q\right)$, with $V=\mathbb{C}\left(K\times C\right)$
(remember the notations $K=\ker\partial$, $I=\mathrm{im}\,\partial$
and $C=\mathrm{coker}\,\partial=\FB{}/I$ from Eq.\ref{eq:exseq})
and $P\left(m\right)\phi:\left(x,Iy\right)\mapsto\delta\left(m,x^{y}\right)\phi\left(x,Iy\right)$,
$Q\left(g\right)\phi:\left(x,Iy\right)\mapsto\phi\left(x,Iyg\right)$
for $\phi\in V$, is an object of $\xrep$, that we'll call the vacuum
object.
\end{example}
Note that the above objects, which exist for any finite crossed module
$\mathcal{X}$, need not be distinct. For example, in the category
$\mathscr{M}\left(\mathcal{R}G\right)$ (see Example \ref{ex:RG})
one has $\vac=\mathbf{\underline{R}}$, while in $\mathscr{M}\left(\mathcal{D}G\right)$
one has $\vac=\idr$.

Given an object $\left(V,P,Q\right)$ of $\xrep$, a linear subspace
$W<V$ is invariant if $P(m)W\subset W$ and $Q(g)W\subset W$ for
all $m\in\FA{}$ and $g\in\FB{}$. An object $\left(V,P,Q\right)$
is reducible if it has a nontrivial invariant subspace, otherwise
it is irreducible. For a finite crossed module $\mathcal{X}$ there
are only finitely many isomorphism classes of irreducible objects
in $\xrep$, which follows from the following generalization of Burnside's
classical theorem \cite{Isaacs,Serre}:

\begin{equation}
\sum_{p\in\irr}d_{p}^{2}=\left|\FB{}\right|\left|\FA{}\right|\,\,,\label{eq:Burnside}\end{equation}
where we denote by $\irr$ the set of (isomorphism classes of) irreducible
objects of $\xrep$, and $d_{p}$ denotes the dimension of the irreducible
$p\in\irr$.

The notion of direct sum of objects of $\xrep$ is the obvious one:
\begin{equation}
\left(V_{1},P_{1},Q_{1}\right)\oplus\left(V_{2},P_{2},Q_{2}\right)=\left(V_{1}\oplus V_{2},P_{1}\oplus P_{2},Q_{1}\oplus Q_{2}\right)\,\,.\label{eq:dirsum}\end{equation}
The analogue of Maschke's theorem states that, for a finite crossed
module $\mathcal{X}$, any object of $\xrep$ decomposes uniquely
(up to ordering) into a direct sum of irreducible objects.

The tensor product of the objects $\left(V_{1},P_{1},Q_{1}\right)$
and $\left(V_{2},P_{2},Q_{2}\right)$ is the triple $\left(V_{1}\otimes V_{2},P_{12},Q_{12}\right)$,
where $P_{12}:m\mapsto\sum_{n\in\FA{}}P_{1}\left(n\right)\otimes P_{2}\left(n^{-1}m\right)$
and $Q_{12}:g\mapsto Q_{1}\left(g\right)\otimes Q_{2}\left(g\right)$.
The category $\xrep$ may be shown to be a monoidal tensor category,
which in general fails to be symmetric, but it is always braided,
the braiding being provided by the map \begin{eqnarray*}
R_{12}: & V_{1}\otimes V_{2}\rightarrow & V_{2}\otimes V_{1}\\
 & v_{1}\otimes v_{2}\mapsto & \sum_{m\in\FA{}}Q_{2}\left(\partial m\right)v_{2}\otimes P_{1}\left(m\right)v_{1}\end{eqnarray*}

At this point it is worthwhile to take a look the category $\xrep$
for the two canonical examples of crossed modules considered in Section
\ref{sec:Crossed-modules}, namely $\mathcal{R}G$ and $\mathcal{D}G$
for a finite group $G$. In the first case, since $\FA{}=\triv$,
the map $P:\FA{}\rightarrow\mathrm{End}\left(V\right)$ is trivial:
$P\left(m\right)=\delta\left(m,1\right)\id$, while the map $Q:\FB{}\rightarrow\mathrm{Aut}\left(V\right)$
provides a representation of the finite group $\FB{}=G$. Thus, for
$\mathcal{X}=\mathcal{R}G$ the category $\xrep$ is nothing but the
category of representations of the finite group $G$. On the other
hand, for $\mathcal{X}=\mathcal{D}G$ the map $P$ is no longer trivial,
and a little thought reveals that in this case $\xrep$ is just the
module category of the (Drinfeld) double of the finite group $G$
\cite{DPR,Hopf1,Hopf2}. It is known that this last tensor category
is modular, and describes the properties of the so-called holomorphic
$G$-orbifold models \cite{DV3}. So, from this point of view, the
category $\xrep$ may be viewed as a common generalization of the
module categories of a finite group and of its double.

\section{Characters\label{sec:Characters}}

The notion of group characters is an extremely powerful tool in the
study of group representations \cite{Isaacs}. Not only do characters
distinguish inequivalent representations, but they prove invaluable
in actual computations, e.g. the decomposition into irreducibles,
the computation of tensor products, etc. As it turns out, a close
analogue of group characters exists for the (isomorphism classes of)
objects of $\xrep$. Namely, the character of an object $\left(V,P,Q\right)$
of $\xrep$ is the complex valued function $\psi:\FA{}\times\FB{}\rightarrow\mathbb{C}$
given by\begin{equation}
\psi\left(m,g\right)=\tr_{V}\left(P\left(m\right)Q\left(g\right)\right)\,\,.\label{chardef}\end{equation}

Clearly, characters of isomorphic objects are equal, and it follows
from the orthogonality relations to be presented a bit later that
characters distinguish inequivalent objects of $\xrep$. The character
$\psi$ of an object of $\xrep$ is a class function of the crossed
module $\mathcal{X}$, i.e. a complex valued function $\psi:\FA{}\times\FB{}\rightarrow\mathbb{C}$
that satisfies

\begin{enumerate}
\item $\psi\left(m,g\right)=0$ unless $m^{g}=m$, for $m\in\FA{}$ and
$g\in\FB{}$;
\item $\psi\left(m^{h},h^{-1}gh\right)=\psi\left(m,g\right)$ for all $m\in\FA{}$
and $g,h\in\FB{}$.
\end{enumerate}
The set of class functions of a finite crossed module $\mathcal{X}$
form a finite dimensional linear space $\cl$, which carries the natural
scalar product \begin{equation}
\left\langle \psi_{1},\psi_{2}\right\rangle =\frac{1}{\left|\FB{}\right|}\sum_{m\in\FA{},g\in\FB{}}\overline{\psi_{1}\left(m,g\right)}\psi_{2}\left(m,g\right)\,\,,\label{clprod}\end{equation}
where $\psi_{1},\psi_{2}\in\cl$, and the bar denotes complex conjugation.

Characters behave well under direct sums and tensor products: the
character of a direct sum is just the (pointwise) sum of the characters
of the summands, while the character of a tensor product is given
by the formula\begin{equation}
\psi_{A\otimes B}\left(m,g\right)=\sum_{n\in\FA{}}\psi_{A}\left(n,g\right)\psi_{B}\left(n^{-1}m,g\right)\,\,,\label{eq:tenschar}\end{equation}
if $\psi_{A},\psi_{B}$ are the characters of the factors.

Irreducible characters, i.e. the characters of the irreducible objects
of $\xrep$, play a distinguished role, since any character may be
written (uniquely) as a linear combination of irreducible ones with
non-negative integer coefficients. The basic result about irreducible
characters is the following analogue of the generalized orthogonality
relations for group characters \cite{Isaacs,Serre}:\begin{equation}
\frac{1}{\left|\FB{}\right|}\sum_{h\in\FB{}}\psi_{p}\left(m,h\right)\psi_{q}\left(m,h^{-1}g\right)=\frac{1}{d_{p}}\delta_{pq}\psi_{p}\left(m,g\right)\label{genorth}\end{equation}
for $p,q\in\irr$, where \begin{equation}
d_{p}=\sum_{m\in\FA{}}\psi_{p}\left(m,1\right)\label{eq:dim}\end{equation}
denotes the dimension of the irreducible $p$. From this one can deduce
at once that the characters of the irreducible representations form
an orthonormal basis in the space $\cl$ of class functions, and that
they also satisfy the second orthogonality relations \begin{equation}
\sum_{p\in\irr}\psi_{p}\left(m,g\right)\psi_{p}\left(n,h\right)=\sum_{z\in\FB{}}\delta\left(n,m^{z}\right)\delta\left(h^{-1},g^{z}\right)\,\,\,.\label{secorth}\end{equation}
Note that the irreducible characters $\psi_{p}$ may be computed explicitly
for any finite crossed module $\mathcal{X}$, e.g. one has $\psi_{\idr}\left(m,g\right)=\delta\left(m,1\right)$
for the identity object $\idr$ of $\xrep$ (cf. Example \ref{ex:identity}).

Using the orthogonality relations, one may express the fusion rule
coefficient $N_{pq}^{r}$, i.e. the multiplicity of the irreducible
$r\in\irr$ in the tensor product of the irreducibles $p$ and $q$,
through the formula \begin{equation}
N_{pq}^{r}=\frac{1}{\left|\FB{}\right|}\sum_{m,n\in\FA{}}\sum_{g\in\FB{}}\psi_{p}\left(m,g\right)\psi_{q}\left(n,g\right)\overline{\psi_{r}\left(mn,g\right)}\,\,.\label{eq:fusionrule}\end{equation}

To each irreducible $p\in\irr$ one may associate the complex number
\begin{equation}
\omega_{p}=\frac{1}{d_{p}}\sum_{m\in\FA{}}\psi_{p}\left(m,\partial m\right)\,\,,\label{omdef}\end{equation}
(remember that $d_{p}$ denotes the dimension of the irreducible $p$),
which turns out to be a root of unity (of order dividing the exponent
of $I=\mathrm{im}\,\partial$), and one may show that \begin{equation}
\psi_{p}\left(m,g\partial m\right)=\omega_{p}\psi_{p}\left(m,g\right)\,\,,\label{omeq}\end{equation}
 for all $m\in\FA{},g\in\FB{}$. Combined with the orthogonality relations
Eq.(\ref{genorth}), this leads to (remember that $K=\ker\partial$)\begin{equation}
\sum_{p\in\irr}d_{p}^{2}\omega_{p}^{-1}=\left|\FB{}\right|\left|K\right|\,\,,\label{eq:dimsum}\end{equation}
to be compared with Eq.(\ref{eq:Burnside}).

To conclude, let's just note that the close analogy with ordinary
group characters goes much further, e.g. one may introduce the Frobenius-Schur
indicator \begin{equation}
\nu_{p}=\frac{1}{\left|\FB{}\right|}\sum_{m\in\FA{},g\in\FB{}}\delta\left(m^{g},m^{-1}\right)\psi_{p}\left(m,g^{2}\right)\,\,.\label{FSdef}\end{equation}
of the irreducible character $\psi_{p}$, and show that $\nu_{p}$
may take only the values $0$ and $\pm1$, in perfect parallel with
the classical case \cite{Isaacs}. Of course, this is related to the
fact that ordinary characters of the finite group $G$ are nothing
but the characters of the crossed module $\mathcal{R}G$ of Example
\ref{ex:RG}.

\section{The $S$ matrix and the structure of the vacuum\label{sec:premodular}}

Up to now, we have seen the close parallel between the structure of
the category $\xrep$ and the representation category of a finite
group. We now turn to describe the premodular structure, related to
the existence of the so-called $S$ matrix. This is a square matrix,
with rows and columns labeled by the irreducibles of $\xrep$, and
with matrix elements \begin{equation}
S_{pq}=\frac{1}{\left|\mathcal{X}\right|}\sum_{m,n\in\FA{}}\overline{\psi_{p}\left(m,\partial n\right)\psi_{q}\left(n,\partial m\right)}\,\,\,\label{sdef}\end{equation}
for $p,q\in\irr$, where $\left|\mathcal{X}\right|=\left|\FA{}\right|\left|C\right|=\left|K\right|\left|\FB{}\right|$
(remember Eq.(\ref{eq:exseq})). This matrix is obviously symmetric,
and a simple computation shows that \begin{equation}
S_{\idr p}=\frac{d_{p}}{\left|\mathcal{X}\right|}>0\,\,,\label{eq:qdim}\end{equation}
where $\idr$ denotes the identity object of $\xrep$ (cf. Example
\ref{ex:identity}).

A most important feature of the above $S$ matrix is its relation
to the fusion rule coefficients $N_{pq}^{r}$ appearing in Eq.(\ref{eq:fusionrule}),
for one may show that \begin{equation}
\sum_{r\in\irr}N_{pq}^{r}S_{rs}=\frac{S_{ps}S_{qs}}{S_{\idr s}}\,\,\label{Ver2}\end{equation}
holds, which is an avatar of Verlinde's celebrated formula \cite{Verlinde_formula}.
A closely related result states that \begin{equation}
\sum_{r\in\irr}N_{pq}^{r}\omega_{r}^{-1}S_{\idr r}=\omega_{p}^{-1}\omega_{q}^{-1}S_{pq}\,\,,\label{eq:Ver1}\end{equation}
where the roots of unity $\omega_{p}$ are given by Eq.(\ref{omdef}).
But this is not the end of the story since, upon introducing the diagonal
matrix $T_{pq}=\omega_{p}\delta_{pq}$, one may show that \begin{equation}
STS=T^{-1}ST^{-1}\,\,\,.\label{modrel}\end{equation}

Should $S$ satisfy the relation $S^{4}=1$, Eq.(\ref{modrel}) would
mean that the matrices $S$ and $T$ give a finite dimensional representation
of the modular group $\mathrm{SL}_{2}\left(\mathbb{Z}\right)$, which
conforms with Verlinde's theorem \cite{Verlinde_formula,Moore-Seiberg},
i.e.

\begin{enumerate}
\item $T$ is diagonal and of finite order;
\item $S$ is symmetric;
\item Verlinde's formula Eq.(\ref{Ver2}) holds.
\end{enumerate}
Should this be the case, $\xrep$ would be a Modular Tensor Category.
As it turns out, in general this is not the case, because the matrix
$S$ of Eq.(\ref{sdef}) does only satisfy the weaker property \begin{equation}
S^{8}=S^{4}\,\,.\label{eq:S4eq}\end{equation}
This means that $S$ is not necessarily invertible: it might have
a nontrivial kernel. This is the extent to which $\xrep$ fails to
be modular in general.

The lack of invertibility of $S$ is related to the reducibility of
the vacuum object $\vac$ (cf. Example \ref{ex:vacuum}). Denoting
by $\mu_{p}$ the multiplicity of the irreducible $p$ in $\vac$,
and by $D=\left|C\right|\left|K\right|$ the dimension of $\vac$,
one may show that \begin{equation}
\mu_{p}=D\left[S^{2}\right]_{\idr p\,\,,}\label{eq:muS2}\end{equation}
and that $\mu_{p}>0$ if and only if there exists an $\alpha$ such
that \begin{equation}
S_{pq}=\alpha S_{\idr q}\,\,\mathrm{for}\,\,\mathrm{all}\,\, q\in\irr\,,\label{eq:deftransparent}\end{equation}
 in which case $\alpha=\mu_{p}=d_{p}$ and $\omega_{p}=1$. In other
words, the irreducible objects of $\xrep$ that satisfy Eq.(\ref{eq:deftransparent})
for some constant $\alpha$ are precisely the irreducible constituents
of the vacuum $\vac$. The invertibility of $S$ requires that the
only such object is the identity $\idr$, and this condition may be
shown to be equivalent to the bijectivity of the boundary map $\partial$,
which in turn is equivalent to $\mathcal{X}$ being isomorphic to
$\mathcal{D}G$ for some finite group $G$. Note also that for $\mathcal{X}=\mathcal{R}G$
every irreducible of $\xrep$ satisfies Eq.(\ref{eq:deftransparent}),
since in this case $\vac=\mathbf{\underline{R}}$.

Finally, we note that while $\xrep$ fails to be modular in case $\partial$
is not bijective, it can nevertheless be turned into an MTC! Indeed,
according to the modularizability criterion of Bruguieres \cite{Bruguieres},
one can associate a well-defined MTC (unique up to isomorphism) to
any premodular category in which Eq.(\ref{eq:deftransparent}) implies
$\omega_{p}=1$ and $\alpha=d_{p}$. But we won't pursue this line
any further in the present note, and leave the construction of the
corresponding MTC to some future work.

\section{Discussion }

As we have sketched in the previous sections, to any finite crossed
module $\mathcal{X}$ one may associate a premodular category $\xrep$.
In special instances this construction gives back the module category
of a finite group or that of its (Drinfeld) double, but in general
one gets new premodular categories, which are very close to being
modular: they satisfy the modularizability criterion of \cite{Bruguieres},
i.e. they can be turned into a Modular Tensor Category. This opens
the way to the construction of a huge number of Modular Tensor Categories
starting from (relatively) simple algebraic structures.

As stressed before, the category $\xrep$ may be viewed as a generalization
of the module category of the double of a finite group $G$, which
describes the properties of holomorphic $G$-orbifolds \cite{DPR,Hopf1,Hopf2}.
This leads to the speculation that for a general crossed module $\mathcal{X}$
the category $\xrep$, or more precisely its modularisation, should
describe the properties of some 'generalized' holomorphic orbifold
related to $\mathcal{X}$. To find out whether this vague idea may
be made to work seems to be a rewarding task.

\bigskip
\emph{Acknowledgments}: This work was supported by research grants
OTKA T047041, T037674, T043582, TS044839, the J\'anos Bolyai
Research Scholarship of the Hungarian Academy of Sciences, and
EC Marie Curie RTN, MRTN-CT-2004-512194.

\bibliographystyle{amsplain}

\end{document}